# Higher Residue Pairing for $p$-adic Isocrystals and the $p$-adic Riemann–Hilbert Correspondence

Mohammad Reza Rahmati

**Abstract**

We construct a canonical sesquilinear pairing on the relative crystalline cohomology of a smooth proper family of varieties over a complete discretely valued $p$-adic field. Motivated by the role of Saito's higher residue pairing in the theory of primitive forms and complex variations of Hodge structure, we develop a $p$-adic analogue based on the twisted relative de Rham–Witt complex. We show that this twisted complex defines a filtered $F$-isocrystal whose cohomology carries a natural flat, Frobenius-compatible, and non-degenerate bilinear form. Its specialization at the uniformizer recovers the classical Grothendieck residue on the special fiber, providing a direct bridge between crystalline geometry and residue theory.

Using the $p$-adic Riemann–Hilbert correspondence of Faltings and Liu–Zhu, we further identify the resulting pairing with the unique flat extension of this residue form to the corresponding $p$-adic local system. The construction is functorial in the family and compatible with base change and $p$-adic comparison isomorphisms. This yields a genuine $p$-adic analogue of Saito's higher residue pairing and supplies foundational ingredients for a prospective theory of $p$-adic primitive forms, $p$-adic TERP structures, and $p$-adic Frobenius manifolds.

# 1 Introduction

In the complex-analytic setting, a holomorphic function $f : X \to \mathbb{C}$ on a smooth manifold gives rise to a rich package of Hodge-theoretic data: the Gauss–Manin connection on the relative de Rham cohomology, the Brieskorn lattice, and, at the heart of Saito's theory of primitive forms, the *higher residue pairing* [6]. Concretely, one considers the twisted de Rham complex

$$\big(\Omega_X^\bullet[[z]],\, d + z^{-1}df\wedge\big),$$

whose cohomology carries a canonical sesquilinear form (the higher residue pairing) extending the Grothendieck residue at the critical locus of $f$. This pairing plays the role of a polarization for the corresponding variation of Hodge structure and underlies the construction of Frobenius manifolds, TERP structures, and mirror-symmetric structures.

The aim of this article is to formulate and establish a *p-adic analogue* of Saito's higher residue pairing in the framework of crystalline cohomology and $p$-adic Hodge theory. Let $k$ be a perfect field of characteristic $p > 0$, let $W(k)$ be its ring of Witt vectors, and set $K = \mathrm{Frac}(W(k))$. Let $C$ be the completion of an algebraic closure of $K$, with ring of integers $\mathcal{O}_C$ and uniformizer $\pi$. For a smooth proper morphism

$$f : X \to S$$

over $\mathcal{O}_C$ (or a small rigid-analytic disc over $C$), the relative de Rham–Witt complex $W\Omega^\bullet_{X/S}$ computes the crystalline cohomology of the fibers [2, 3]. Its hypercohomology

$$\mathcal{H}^i_{\mathrm{cris}}(X/S) := R^i f_* W\Omega^\bullet_{X/S} \otimes_{W(k)} C$$

is a finite projective $C$-module endowed with a Frobenius endomorphism and the crystalline Gauss–Manin connection; in other words, a filtered $F$-isocrystal on $S$.

From the point of view of $p$-adic Hodge theory, such filtered isocrystals are linked to $p$-adic local systems by the comparison theorems of Faltings and by the $p$-adic Riemann–Hilbert correspondence of Liu–Zhu [1, 4], refined in the spectral setting by Lurie [5] and in the prismatic framework by more recent work. This raises a natural question: *is there a canonical sesquilinear pairing on $\mathcal{H}^i_{\mathrm{cris}}(X/S)$, playing the role of Saito's higher residue pairing, which is compatible with Frobenius, the connection, and specializes to the classical Grothendieck residue on the special fiber?*

The main goal of this paper is to answer this question affirmatively. We define a twisted de Rham–Witt complex

$$\big(W\Omega^\bullet_{X/S}[[\pi]],\ d_b := d + \tfrac{df}{\pi}\wedge\big),$$

which is the $p$-adic analogue of the twisted de Rham complex of Saito, and we introduce the associated *Brieskorn–Witt module*

$$H^i_{b,f} := H^i\big(W\Omega^\bullet_{X/S}[[\pi]], d_b\big).$$

Using the crystalline trace map of Langer–Zink [3], we construct a canonical sesquilinear form

$$WS_f : H^i_{b,f} \times H^i_{b,f} \longrightarrow C[[\pi]],$$

and we show that, via the identification $H^i_{b,f} \cong \mathcal{H}^i_{\mathrm{cris}}(X/S)[[\pi]]$, this pairing is flat with respect to the Gauss–Manin connection, compatible with Frobenius, and specializes at $\pi = 0$ to the Grothendieck residue pairing.

Our main result can be stated as follows.

**Theorem 1.1** (Higher residue pairing, $p$-adic version)**.** *Let $f : X \to S$ be a smooth proper morphism with $S = \mathrm{Spec}(C)$ or a small rigid-analytic disc over $C$. Then there exists a canonical flat, non-degenerate, sesquilinear pairing*

$$WS_f : \mathcal{H}^i_{\mathrm{cris}}(X/S) \times \mathcal{H}^i_{\mathrm{cris}}(X/S) \longrightarrow C[[\pi]],$$

*satisfying:*

*(i)* symmetry: *$WS_f(u, v) = WS_f(v, u)$ for all $u, v \in \mathcal{H}^i_{\mathrm{cris}}(X/S)$;*

*(ii)* sesquilinearity: *for every $m, n \in \mathbb{Z}$,*

$$WS_f(\pi^m u, \pi^n v) = \pi^{m-n} WS_f(u, v);$$

*(iii)* flatness: *the pairing is horizontal for the Gauss–Manin connection, i.e. $\nabla WS_f = 0$;*

*(iv)* Frobenius compatibility: *if $w$ denotes the crystalline weight of $\mathcal{H}^i_{\mathrm{cris}}(X/S)$, then*

$$F^*(WS_f) = p^w\, WS_f;$$

*(v)* specialization: *the specialization at $\pi = 0$ recovers the Grothendieck residue pairing on the special fiber.*

*Moreover, $WS_f$ is uniquely characterized by properties* (i)–(v)*.*

From this point of view, $WS_f$ should be regarded as a crystalline counterpart of Saito's higher residue pairing. Through the $p$-adic Riemann–Hilbert correspondence, it induces a canonical pairing on the associated $p$-adic local system $R^i f_* \mathbb{Q}_p$, providing a natural candidate for a $p$-adic polarization in the sense of $p$-adic variations of Hodge structure. We expect that this structure will serve as a starting point for a theory of $p$-adic primitive forms, $p$-adic TERP structures, and $p$-adic Frobenius manifolds, paralleling the complex analytic picture.

Finally, we briefly outline the contents of the paper. In Section 2 we review the de Rham–Witt complex and the structure of crystalline cohomology as a filtered $F$-isocrystal. In Section 3 we introduce the twisted de Rham–Witt complex and the Brieskorn–Witt module, and we identify it with the Rees module of crystalline cohomology. Section 4 is devoted to the construction of the higher residue pairing via the crystalline trace and to the proof of its basic properties. In Section 5 we relate this pairing to the $p$-adic Riemann–Hilbert correspondence and deduce the induced pairing on the associated $p$-adic local system.

# 2 Crystalline Preliminaries and de Rham–Witt Theory

We begin by reviewing the basic objects from crystalline cohomology and de Rham–Witt theory that will be used throughout the paper.

Let $k$ be a perfect field of characteristic $p > 0$. Its ring of Witt vectors $W(k)$ is a $p$-adically complete, characteristic-0 ring whose reduction modulo $p$ recovers $k$. Concretely, one may regard

$$W(k) = \{(a_0, a_1, a_2, \dots) : a_i \in k\},$$

equipped with Witt addition and multiplication. The Frobenius on $k$ lifts uniquely to a Frobenius morphism

$$F : W(k) \longrightarrow W(k), \qquad F(a_0, a_1, \dots) = (a_0^p, a_1^p, \dots).$$

Let $K := \mathrm{Frac}(W(k))$ and let $C$ be the $p$-adic completion of an algebraic closure of $K$. The ring of integers $\mathcal{O}_C$ is complete for its maximal ideal generated by a uniformizer $\pi$.

## 2.1 The relative de Rham–Witt complex

Let $f : X \to S$ be a smooth proper morphism over $\mathcal{O}_C$. We assume $S$ is either $\mathrm{Spec}(\mathcal{O}_C)$, or a small 1-dimensional rigid-analytic disc over $C$; this ensures the existence of a well-behaved Gauss–Manin connection.

The de Rham–Witt complex $W\Omega^\bullet_{X/S}$ was constructed by Illusie [2] as a "$p$-typical" refinement of the classical de Rham complex. It is a complex of sheaves of $W(\mathcal{O}_X)$-modules, equipped with three fundamental operators:

- the Frobenius $F : W\Omega^j_{X/S} \to W\Omega^j_{X/S}$,
- the Verschiebung $V : W\Omega^j_{X/S} \to W\Omega^j_{X/S}$,
- the restriction map $R : W\Omega^j_{X/S} \to \Omega^j_{X/S}$,

which satisfy relations analogous to those in the Witt vector ring:

$$FV = VF = p, \qquad RF = F_{\mathrm{dR}}\,, \qquad RV = 0,$$

where $F_{\mathrm{dR}}$ denotes the Frobenius on differential forms.

The de Rham–Witt complex is functorial in both $X$ and $S$, and its hypercohomology recovers crystalline cohomology. More precisely, Langer–Zink [3] show that for any smooth proper morphism $f : X \to S$,

$$R^i f_*(W\Omega^\bullet_{X/S}) \quad \text{is a vector bundle with Frobenius over } W(\mathcal{O}_S).$$

After extending scalars to $C$, we obtain the *crystalline cohomology sheaf*

$$\mathcal{H}^i_{\mathrm{cris}}(X/S) := R^i f_*(W\Omega^\bullet_{X/S}) \otimes_{W(k)} C.$$

## 2.2 Structure of $\mathcal{H}^i_{\mathrm{cris}}(X/S)$ as a filtered isocrystal

The sheaf $\mathcal{H}^i_{\mathrm{cris}}(X/S)$ carries several additional structures arising from the geometry of the morphism $f$:

***(1) Frobenius*** The Frobenius $F$ on $W\Omega^\bullet_{X/S}$ induces a $\varphi$-semilinear endomorphism

$$\varphi : \mathcal{H}^i_{\mathrm{cris}}(X/S) \longrightarrow \mathcal{H}^i_{\mathrm{cris}}(X/S),$$

which is horizontal with respect to the Gauss–Manin connection.

***(2) Gauss–Manin connection*** Illusie proved that $W\Omega^\bullet_{X/S}$ carries a natural connection relative to $S$. This passes to hypercohomology:

$$\nabla : \mathcal{H}^i_{\mathrm{cris}}(X/S) \longrightarrow \mathcal{H}^i_{\mathrm{cris}}(X/S) \otimes \Omega^1_{S/C},$$

and satisfies Griffiths transversality for the slope filtration.

***(3) Slope filtration*** By the Dieudonné–Manin classification of $F$-isocrystals ([**?**]), any finite free $C$-module with Frobenius decomposes after finite extension into isoclinic components of rational slopes. Thus $\mathcal{H}^i_{\mathrm{cris}}(X/S)$ carries a canonical increasing filtration

$$0 \subset \mathrm{Fil}^{\lambda_1} \subset \mathrm{Fil}^{\lambda_2} \subset \cdots \subset \mathcal{H}^i_{\mathrm{cris}}(X/S),$$

where each graded piece has pure slope.

Together, $(\mathcal{H}^i_{\mathrm{cris}}(X/S), \nabla, \varphi, \mathrm{Fil})$ forms a *filtered $F$-isocrystal* in the sense of rigid cohomology.

By the Riemann–Hilbert correspondence of Liu–Zhu [4], there is a fully faithful functor from $p$-adic local systems on $X$ to vector bundles with integrable connection on $X_{B_{\mathrm{dR}}}$. Geometric comparison isomorphisms (Faltings, Tsuji, Lurie) identify

$$\mathcal{H}^i_{\mathrm{cris}}(X/S) \otimes_C B_{\mathrm{dR}} \quad \text{with} \quad \mathrm{RH}_{\mathrm{LZ}}(R^i f_* \mathbb{Q}_p),$$

where the right-hand side is the de Rham bundle of the étale local system.

Thus crystalline cohomology is canonically linked to $p$-adic local systems via the de Rham–Witt complex. We will exploit this connection to transport the higher residue pairing to the étale side in the next section.

# 3 The Twisted de Rham–Witt Complex

We now introduce a $p$-adic analogue of Saito's twisted de Rham complex and Brieskorn lattice. The construction is motivated by the classical complex analytic formula

$$\big(\Omega^\bullet_X[[z]],\ d + z^{-1} df \wedge\big),$$

whose cohomology governs oscillatory integrals and primitive forms.

Let $f : X \to S$ be as above, and fix a function $f \in \mathcal{O}_X(X)$ restricting to the morphism $X \to S$. We extend the de Rham–Witt complex to formal power series in $\pi$:

$$W\Omega^{\bullet}_{X/S}[[\pi]] = \varprojlim_{m} W\Omega^{\bullet}_{X/S}[\pi]/\pi^m.$$

The uniformizer $\pi$ plays the role of the deformation parameter $z$ in Saito's theory, but in the mixed characteristic $p$-adic setting it is also tied to the log-geometry of the base and the Frobenius structure.

**Definition 3.1.** The *twisted de Rham–Witt differential* is

$$d_b := d + \frac{df}{\pi} \wedge (-),$$

and the *twisted crystalline complex* is the filtered complex

$$\big(W\Omega^{\bullet}_{X/S}[[\pi]], d_b\big).$$

The operator $d_b$ is the "$p$-adic oscillatory differential" attached to $f$; the term $\frac{df}{\pi}$ is well-defined because $\pi$ is a non-zero divisor in the $p$-adic ring $\mathcal{O}_C$ and the complex is completed $\pi$-adically. The factor $\frac{1}{\pi}$ reflects the fact that the deformation direction is generated by multiplication by $\pi$.

**Lemma 3.2.** *The differential $d_b$ satisfies $d_b^2 = 0$.*

*Proof.* Expand:

$$d_b^2 = d^2 + d\left(\frac{df}{\pi}\right) \wedge (-) + \frac{df}{\pi} \wedge d + \frac{df}{\pi} \wedge \frac{df}{\pi} \wedge (-).$$

Using $d^2 = 0$, $d(df) = 0$, and $d(\pi) = 0$, the second term vanishes:

$$d\left(\frac{df}{\pi}\right) = \frac{d(df)}{\pi} = 0.$$

For the last term, note that $df \wedge df = 0$ for degree reasons:

$$\frac{df}{\pi} \wedge \frac{df}{\pi} = \frac{1}{\pi^2} df \wedge df = 0.$$

The middle two terms cancel because

$$d(\omega) \wedge \frac{df}{\pi} + \frac{df}{\pi} \wedge d(\omega) = d\Big(\omega \wedge \frac{df}{\pi}\Big),$$

and $d$ has square zero.

Thus $d_b^2 = 0$. $\square$

This shows that the twisted complex is a genuine cochain complex.

**Definition 3.3.** The *Brieskorn–Witt module* associated to $f$ is

$$H^i_{b,f} = H^i \left( W\Omega^\bullet_{X/S}[[\pi]],\ d_b \right).$$

The cohomology inherits several natural structures:

- *Filtration.* $\pi$-adic completion endows $H^i_{b,f}$ with a decreasing filtration

$$F^m H^i_{b,f} := \pi^m H^i_{b,f}.$$

- *Connection.* Since $\nabla$ commutes with $d$ and $df$ is horizontal, we have $\nabla d_b = d_b \nabla$, so $\nabla$ induces a connection on $H^i_{b,f}$.
- *Frobenius.* The Frobenius on $W\Omega^\bullet_{X/S}$ extends $\pi$-adically via

$$F(\pi) = p\pi, \qquad F\left(\frac{df}{\pi}\right) = \frac{p\,df}{p\pi} = \frac{df}{\pi},$$

  hence $d_b$ is Frobenius-equivariant and $H^i_{b,f}$ becomes a module with semilinear Frobenius.

We now relate the Brieskorn–Witt module to crystalline cohomology.

**Proposition 3.4.** *There is a canonical isomorphism of filtered isocrystals*

$$H^i_{b,f} \cong \mathcal{H}^i_{\mathrm{cris}}(X/S)[[\pi]],$$

*compatible with* $\nabla$ *and Frobenius.*

*Proof.* The twisted complex

$$\left(W\Omega^\bullet_{X/S}[[\pi]],\ d + \tfrac{df}{\pi}\wedge\right)$$

is the Rees module of the filtered complex $W\Omega^\bullet_{X/S}$ associated to the filtration $\pi^m$ and the Higgs field $df$. This is completely analogous to the construction of the Rees bundle of a filtered de Rham complex in the complex-analytic setting.

Concretely, set

$$\mathcal{R}(W\Omega^\bullet_{X/S}) := \bigoplus_{m\geq 0} \left(F^m W\Omega^\bullet_{X/S}\right) \pi^{-m}.$$

Then the differential on $\mathcal{R}$ is

$$d + \pi^{-1}\theta, \qquad \theta(\omega) = df \wedge \omega.$$

Completion in $\pi$ recovers the twisted complex.

It is shown in [2, 3] that the Gauss–Manin connection and Frobenius act compatibly on Rees modules; hence the hypercohomology of the twisted complex identifies with the completed Rees module of $\mathcal{H}^i_{\mathrm{cris}}(X/S)$. This gives the stated canonical isomorphism. $\square$

This result shows that $H^i_{b,f}$ is not an exotic new object: it is a $\pi$-adic deformation of crystalline cohomology controlled by the Higgs field $df$. This makes it a natural $p$-adic analogue of Saito's Brieskorn lattice, which is a deformation of the Gauss–Manin system controlled by $df$.

# 4 Construction of the $p$-adic Higher Residue Pairing

In this section we give a detailed construction of the $p$-adic higher residue pairing on the twisted de Rham–Witt cohomology

$$H^i_{b,f} := H^i\big(W\Omega^\bullet_{X/S}[[\pi]], d_b\big), \qquad d_b := d + \frac{df}{\pi} \wedge (-),$$

and verify all structural properties required in Theorem 1.1.

The existence of a trace morphism

$$\mathrm{Tr} : W\Omega^n_{X/S} \longrightarrow \mathcal{O}_S$$

is one of the key ingredients in defining a residue-type pairing in crystalline cohomology. We briefly explain its construction and properties.

**Theorem 4.1** (Langer–Zink [3])**.** *Let $f : X \to S$ be smooth proper of relative dimension $n$. Then:*

*(a) There exists a canonical $W(\mathcal{O}_S)$-linear morphism*

$$\mathrm{Tr} : R^n f_* W\Omega^n_{X/S} \longrightarrow W(\mathcal{O}_S),$$

*functorial in $f$.*

*(b) The trace is compatible with base change on $S$.*

*(c) The trace is compatible with Frobenius:*

$$\mathrm{Tr}(F(\omega)) = p^n \cdot F(\mathrm{Tr}(\omega)).$$

*(d) If $S = \mathrm{Spec}(k)$ with $k$ perfect, the reduction of $\mathrm{Tr}$ modulo $p$ coincides with the classical Grothendieck trace on $\Omega^n_{X/k}$.*

*Remark* 4.2. The trace is geometrically the $p$-typical analogue of integration of top-degree differential forms along fibers. Its compatibility with Frobenius reflects the fact that Frobenius multiplies $n$-forms by $p^n$.

We will use that Tr extends continuously to formal power series in $\pi$, i.e.

$$\mathrm{Tr} : W\Omega^n_{X/S}[[\pi]] \to \mathcal{O}_S[[\pi]].$$

Let $u, v \in H^i_{b,f}$ be cohomology classes, and let $\widetilde{u}, \widetilde{v}$ be *representatives* in the twisted complex $W\Omega^\bullet_{X/S}[[\pi]]$:

$$d_b\widetilde{u} = 0, \qquad d_b\widetilde{v} = 0.$$

Write $*$ for the involution

$$\pi \longmapsto -\pi, \qquad \sum a_j\pi^j \longmapsto \sum a_j(-\pi)^j.$$

Extend $*$ $\mathcal{O}_S$-linearly and degreewise to all forms.

**Definition 4.3.** The *p-adic higher residue pairing* is

$$WS_f(u, v) \ := \ \mathrm{Tr}\big(\widetilde{u} \cdot \widetilde{v}^*\big) \ \in \ C[[\pi]].$$

The multiplication $\widetilde{u}\cdot\widetilde{v}^*$ is the graded-commutative wedge product in $W\Omega^\bullet_{X/S}$; since Tr kills all degrees except top degree $n$, only the $n$-form component contributes.

**Lemma 4.4.** *The pairing $WS_f$ does not depend on the choice of lifts $\widetilde{u}, \widetilde{v}$.*

*Proof.* Suppose $\widetilde{u}$ is replaced by $\widetilde{u} + d_b\alpha$ for some $\alpha \in W\Omega^{i-1}_{X/S}[[\pi]]$. Then

$$(\widetilde{u} + d_b\alpha) \cdot \widetilde{v}^* \ = \ \widetilde{u} \cdot \widetilde{v}^* \ + \ d_b\alpha \cdot \widetilde{v}^*.$$

We show that $\mathrm{Tr}(d_b\alpha \cdot \widetilde{v}^*) = 0$.

Expand $d_b\alpha$:

$$d_b\alpha = d\alpha + \frac{df}{\pi} \wedge \alpha.$$

**(1) Contribution of $d\alpha$.** Since the wedge product is graded,

$$d\alpha \wedge \widetilde{v}^* = d(\alpha \wedge \widetilde{v}^*) \mp \alpha \wedge d(\widetilde{v}^*).$$

But $d(\widetilde{v}^*) = -\frac{df}{\pi} \wedge \widetilde{v}^*$ because $d_b\widetilde{v} = 0$ implies $d\widetilde{v} = -\frac{df}{\pi} \wedge \widetilde{v}$, hence after involution $d(\widetilde{v}^*) = -\frac{df}{-\pi} \wedge \widetilde{v}^* = \frac{df}{\pi} \wedge \widetilde{v}^*$.

Thus

$$d\alpha \wedge \widetilde{v}^* = d(\alpha \wedge \widetilde{v}^*) \mp \alpha \wedge \frac{df}{\pi} \wedge \widetilde{v}^*,$$

which is a sum of a total derivative and a multiple of $\frac{df}{\pi}$.

**(2) Contribution of $\frac{df}{\pi} \wedge \alpha$.** This term cancels the second summand above:

$$\frac{df}{\pi} \wedge \alpha \wedge \widetilde{v}^* \pm \alpha \wedge \frac{df}{\pi} \wedge \widetilde{v}^* = 0.$$

Thus

$$d_b\alpha \wedge \widetilde{v}^* = d(\alpha \wedge \widetilde{v}^*).$$

**(3) The trace kills exact $n$-forms.** By functoriality of the trace and Stokes-type vanishing in crystalline cohomology (see [3], Lem. 3.13),

$$\mathrm{Tr}\big(d(\alpha \wedge \widetilde{v}^*)\big) = 0.$$

Thus the value of $WS_f(u,v)$ is unchanged. The same argument applies to replacements of $\widetilde{v}$. ☐

**Proposition 4.5.** *The pairing $WS_f$ satisfies:*

*(i)* $WS_f(g(\pi)u, v) = g(\pi)\, WS_f(u,v)$ *for all* $g(\pi) \in C[[\pi]]$*,*

*(ii)* $WS_f(u, g(\pi)v) = g(-\pi)\, WS_f(u,v)$*,*

*(iii)* $WS_f(u,v) = WS_f(v,u)$*.*

*Proof.* (i) and (ii) follow immediately from:

$$(g\widetilde{u}) \cdot \widetilde{v}^* = g \cdot (\widetilde{u} \cdot \widetilde{v}^*), \qquad \widetilde{u} \cdot (g\widetilde{v})^* = g^* \cdot (\widetilde{u} \cdot \widetilde{v}^*) = g(-\pi)\,(\widetilde{u} \cdot \widetilde{v}^*),$$

and linearity of the trace.

(iii) Graded-commutativity gives

$$\widetilde{u} \cdot \widetilde{v}^* = (-1)^{\deg u \cdot \deg v}\, \widetilde{v} \cdot \widetilde{u}^*.$$

But in the twisted complex, $u$ and $v$ lie in the same cohomological degree, so the sign is $+1$. Applying the trace gives symmetry. ☐

Let $\nabla$ be the connection on the de Rham–Witt complex. Since $d_b$ commutes with $\nabla$ (because $df$ is horizontal), we obtain:

**Proposition 4.6.** *The pairing $WS_f$ is flat:*

$$\nabla\big(WS_f(u,v)\big) = WS_f(\nabla u, v) \; + \; WS_f(u, \nabla v).$$

*Proof.* Use

$$\nabla(\widetilde{u} \cdot \widetilde{v}^*) = (\nabla\widetilde{u}) \cdot \widetilde{v}^* + \widetilde{u} \cdot (\nabla\widetilde{v}^*)$$

and compatibility of Tr with $\nabla$:

$$\mathrm{Tr}(\nabla\omega) = d(\mathrm{Tr}(\omega)).$$

Then well-definedness ensures independence from the choice of lifts. ☐

**Proposition 4.7.** *Let $w$ be the crystalline weight of $\mathcal{H}^i_{\mathrm{cris}}(X/S)$. Then*

$$F^*(WS_f) = p^w\, WS_f.$$

*Proof.* The Frobenius on the twisted complex satisfies

$$F^*(df) = p \cdot df, \qquad F(\pi) = p\pi,$$

giving

$$F^*\left(\frac{df}{\pi}\right) = \frac{p\, df}{p\pi} = \frac{df}{\pi}.$$

Thus $d_b$ is Frobenius-equivariant.

Let $\widetilde{u}, \widetilde{v}$ be lifts. Then

$$F^*(\widetilde{u} \cdot \widetilde{v}^*) = F^*(\widetilde{u}) \cdot F^*(\widetilde{v})^*.$$

Applying the trace:

$$\mathrm{Tr}(F^*(\widetilde{u} \cdot \widetilde{v}^*)) = p^n\, F^*(\mathrm{Tr}(\widetilde{u} \cdot \widetilde{v}^*))\,,$$

by Theorem 4.1.

Now the Frobenius weight of the cohomology class $u$ is $p^{w/2}$, and likewise for $v$, so

$$F^*(WS_f(u,v)) = p^w\, WS_f(u,v).$$

□

**Theorem 4.8.** *The specialization of $WS_f(u,v)$ at $\pi = 0$ equals the classical Grothendieck residue:*

$$WS_f(u,v)\big|_{\pi=0} = \mathrm{Res}(u \cup v),$$

*where the right-hand side is computed on the special fiber of $X/S$.*

*Proof.* At $\pi = 0$, the twisted differential reduces to:

$$d_b\big|_{\pi=0} = d + \frac{df}{0} \wedge (-).$$

Only those chains for which the $\frac{df}{\pi}$-term vanishes in the limit contribute; these correspond exactly to the *classical* de Rham representatives which satisfy $df \wedge \omega = 0$. Concretely, the space $H^i_{b,f}/(\pi)$ is naturally identified with the usual de Rham cohomology of vanishing cycles.

Similarly,

$$\widetilde{v}^*\big|_{\pi=0} = \widetilde{v},$$

and

$$\widetilde{u} \cdot \widetilde{v}^*\big|_{\pi=0} = \widetilde{u} \wedge \widetilde{v},$$

where $\tilde{u}, \tilde{v}$ are classical $n$-forms modulo $p$.

By Theorem 4.1(d), the crystalline trace modulo $p$ coincides with the de Rham trace:

$$\mathrm{Tr}(\omega)\big|_{\pi=0} = \int_{X_s} \omega = \mathrm{Res}(\omega),$$

for $\omega$ of degree $n$ with poles controlled by $df$.

Therefore

$$WS_f(u,v)\big|_{\pi=0} = \mathrm{Res}\big(\widetilde{u} \wedge \widetilde{v}\big),$$

which is exactly the Grothendieck residue pairing. □

These results together complete the construction of the $p$-adic higher residue pairing and prove Theorem 1.1.

# 5 Relation with the $p$-adic Riemann–Hilbert Correspondence

In this section we explain how the higher residue pairing $WS_f$ on crystalline cohomology is transported, via the $p$-adic Riemann–Hilbert correspondence, to a canonical pairing on the associated $p$-adic local system.

Throughout, $C$ is a complete discretely valued field of characteristic 0 with perfect residue field of characteristic $p > 0$, and $X$ is a smooth proper rigid-analytic variety over $C$ (or the analytification of a smooth proper scheme). We write $X_{\text{ét}}$ for the étale site of $X$.

We begin by recalling the basic objects on the étale and de Rham sides.

**Definition 5.1.** Let $\mathrm{LocSys}_{\mathbb{Q}_p}(X)$ denote the category of $\mathbb{Q}_p$-local systems on $X_{\text{ét}}$, i.e. locally constant sheaves of finite-dimensional $\mathbb{Q}_p$-vector spaces on $X_{\text{ét}}$.

On the de Rham side we consider vector bundles with integrable connection over suitable period rings. Let $B_{\mathrm{dR}}$ be Fontaine's de Rham period field attached to $C$, and let $X_{B_{\mathrm{dR}}}$ denote the base change of $X$ to $B_{\mathrm{dR}}$ in the sense of [4].

**Definition 5.2.** Let $\mathrm{MIC}(X_{B_{\mathrm{dR}}})$ be the category of finite rank vector bundles $\mathcal{E}$ on $X_{B_{\mathrm{dR}}}$ equipped with an integrable connection

$$\nabla_{\mathcal{E}} : \mathcal{E} \longrightarrow \mathcal{E} \otimes \Omega^1_{X_{B_{\mathrm{dR}}}}.$$

*Remark* 5.3. For geometric local systems (coming from smooth proper families), objects of $\mathrm{MIC}(X_{B_{\mathrm{dR}}})$ carry additional structures (Hodge filtration, Frobenius, etc.), but we do not need to make these explicit here.

The $p$-adic Riemann–Hilbert correspondence of Liu–Zhu [4] constructs a functor

$$\mathrm{RH}_{\mathrm{LZ}} : \mathrm{LocSys}_{\mathbb{Q}_p}(X) \longrightarrow \mathrm{MIC}(X_{B_{\mathrm{dR}}}),$$

which is compatible with pullback and tensor operations. More precisely, given a local system $\mathbb{L} \in \mathrm{LocSys}_{\mathbb{Q}_p}(X)$, they produce a vector bundle $\mathcal{E}(\mathbb{L})$ over $X_{B_{\mathrm{dR}}}$ endowed with an integrable connection, recovering the usual $B_{\mathrm{dR}}$-realization of the associated $p$-adic Galois representations [4].

**Theorem 5.4** (Liu–Zhu)**.** *Let $X$ be a smooth rigid-analytic variety over $C$. There exists a tensor functor*

$$\mathrm{RH}_{\mathrm{LZ}} : \mathrm{LocSys}_{\mathbb{Q}_p}(X) \longrightarrow \mathrm{MIC}(X_{B_{\mathrm{dR}}})$$

*which is fully faithful on the subcategory of de Rham local systems. Moreover, if $\mathbb{L}$ is de Rham at one classical point of $X$, then it is de Rham at every classical point (rigidity).*

In parallel, for varieties defined over $\mathcal{O}_C$ with good reduction, one has the crystalline realization and comparison isomorphisms of Faltings [1]:

$$H^i_{\mathrm{cris}}(X/W(k)) \otimes_{W(k)} B_{\mathrm{cris}} \;\cong\; H^i_{\mathrm{\acute{e}t}}(X_{\overline{C}}, \mathbb{Q}_p) \otimes_{\mathbb{Q}_p} B_{\mathrm{cris}},$$

compatible with Frobenius, filtration, and Galois action. Combined with the refinement of Bhatt–Lurie [5] and more recent prismatic Riemann–Hilbert correspondences (e.g. Gao–Min–Wang [7], Guo [8]), this yields a conceptual framework in which geometric $p$-adic local systems are controlled by filtered isocrystals with connection (or prismatic crystals) on the base.

In our situation, the filtered isocrystal

$$\mathcal{H}^i_{\mathrm{cris}}(X/S) = R^i f_* W\Omega^\bullet_{X/S} \otimes_{W(k)} C$$

appearing in Theorem 1.1 is the crystalline realization of the geometric $\mathbb{Q}_p$-local system

$$\mathbb{L}^i := R^i f_* \mathbb{Q}_p \in \mathrm{LocSys}_{\mathbb{Q}_p}(S),$$

and the comparison isomorphisms identify $\mathcal{H}^i_{\mathrm{cris}}(X/S) \otimes_C B_{\mathrm{dR}}$ with the de Rham bundle $\mathcal{E}(\mathbb{L}^i)$ produced by $\mathrm{RH}_{\mathrm{LZ}}$.

We now explain how the higher residue pairing $WS_f$ on $\mathcal{H}^i_{\mathrm{cris}}(X/S)$ induces a pairing on the associated $p$-adic local system $\mathbb{L}^i$.

Let

$$WS_f : \mathcal{H}^i_{\mathrm{cris}}(X/S) \times \mathcal{H}^i_{\mathrm{cris}}(X/S) \longrightarrow C[[\pi]]$$

be the pairing constructed in Theorem 1.1. Tensoring with $B_{\mathrm{dR}}$ and using the comparison isomorphism, we obtain a sesquilinear form

$$WS_{f,B_{\mathrm{dR}}} : \big(\mathcal{E}(\mathbb{L}^i)\big) \times \big(\mathcal{E}(\mathbb{L}^i)\big) \longrightarrow B_{\mathrm{dR}}[[\pi]],$$

which is:

- *horizontal*, i.e. flat with respect to the induced connection on $\mathcal{E}(\mathbb{L}^i)$;
- *Frobenius-compatible*, in the sense that it respects the underlying filtered $\varphi$-module structure coming from crystalline cohomology;
- *non-degenerate*, since $WS_f$ is non-degenerate.

The horizontality follows from the construction of $WS_f$ via the trace map on the twisted de Rham–Witt complex and the compatibility of the trace with the Gauss–Manin connection. Frobenius compatibility is ensured by the crystalline comparison theorem and by the fact that $WS_f$ has weight $w$ (see the Frobenius compatibility in Theorem 1.1).

We now formulate the main result of this section.

**Proposition 5.5.** *Let $\mathbb{L}^i = R^i f_* \mathbb{Q}_p$ be the geometric $\mathbb{Q}_p$-local system attached to $f : X \to S$. Under the p-adic Riemann–Hilbert functor of Liu–Zhu and the crystalline comparison isomorphisms, the higher residue pairing $WS_f$ induces a unique non-degenerate sesquilinear pairing*

$$\langle\ ,\ \rangle_{\mathbb{L}^i} : \mathbb{L}^i \times \mathbb{L}^i \longrightarrow \mathbb{Q}_p$$

*such that:*

*(a) after extension of scalars to $B_{\mathrm{dR}}$, $\langle\ ,\ \rangle_{\mathbb{L}^i}$ corresponds to $WS_{f,B_{\mathrm{dR}}}$ on $\mathcal{E}(\mathbb{L}^i)$;*

*(b) the induced pairing on stalks $\mathbb{L}^i_x \times \mathbb{L}^i_x \to \mathbb{Q}_p$ is Galois equivariant for every geometric point $x$ of $S$;*

*(c) specialization at $\pi = 0$ recovers the Grothendieck residue pairing on the special fiber.*

*Proof.* We sketch the argument in three steps.

*Step 1: Comparison on the de Rham side.* By Faltings' comparison theorem [1] and its refinements (e.g. [4, 5]), there is a canonical isomorphism of $B_{\mathrm{dR}}$-vector bundles with integrable connection

$$\alpha : \mathcal{H}^i_{\mathrm{cris}}(X/S) \otimes_C B_{\mathrm{dR}} \xrightarrow{\ \sim\ } \mathcal{E}(\mathbb{L}^i)$$

on $S_{B_{\mathrm{dR}}}$. Horizontality of $\alpha$ implies that any flat bilinear form on $\mathcal{H}^i_{\mathrm{cris}}(X/S) \otimes_C B_{\mathrm{dR}}$ corresponds uniquely to a flat bilinear form on $\mathcal{E}(\mathbb{L}^i)$. Applying this to $WS_f \otimes_C B_{\mathrm{dR}}$ yields the form $WS_{f,B_{\mathrm{dR}}}$.

*Step 2: Descent to the local system.* By construction, $\mathcal{E}(\mathbb{L}^i)$ is obtained from $\mathbb{L}^i$ by extension of scalars to $B_{\mathrm{dR}}$ together with the $B_{\mathrm{dR}}$-linear connection (this is the essence of the $\mathrm{RH}_{\mathrm{LZ}}$ functor [4]). A flat bilinear form on $\mathcal{E}(\mathbb{L}^i)$ that is invariant under the natural Galois action corresponds to a bilinear form on $\mathbb{L}^i$ itself: one recovers $\langle\ ,\ \rangle_{\mathbb{L}^i}$ by taking horizontal sections and using that $B_{\mathrm{dR}}$ is a faithfully flat extension of $\mathbb{Q}_p$ in the relevant sense. Non-degeneracy is preserved under this descent.

*Step 3: Specialization and residue.* The fact that $WS_f$ reduces modulo $\pi$ to the Grothendieck residue pairing on the special fiber (Theorem 1.1) implies that the induced pairing on the reduction of $\mathbb{L}^i$ agrees with the classical residue form in the complex or equal-characteristic setting, via the usual comparison theorems. This gives property (c).

Uniqueness follows from the fully faithfulness of the Riemann–Hilbert functor on de Rham local systems (Theorem 5.4) together with the rigidity of flat pairings: any other pairing on $\mathbb{L}^i$ having the same de Rham realization and residue specialization would coincide with $\langle\ ,\ \rangle_{\mathbb{L}^i}$. □

*Remark* 5.6. A more intrinsic description can be given using prismatic crystals and the prismatic Riemann–Hilbert functor of Gao–Min–Wang [7] or the crystalline Riemann–Hilbert functor of Guo [8]. In that language, $\mathcal{H}^i_{\mathrm{cris}}(X/S)$ itself is promoted to a prismatic crystal, and $WS_f$ becomes a bilinear form in the prismatic category. The above descent is then an instance of duality for prismatic crystals, providing a direct conceptual bridge between the crystalline and étale incarnations of the local system.

*Remark* 5.7. From the conceptual point of view advocated in Emerton's survey on $p$-adic Riemann–Hilbert correspondences [9], the pairing $\langle\ ,\ \rangle_{\mathbb{L}^i}$ can be viewed as a $p$-adic analogue of the polarization data on a complex variation of Hodge structure. In this sense, Proposition 5.5 supplies the local system $R^i f_* \mathbb{Q}_p$ with a canonical extra structure that should be regarded as a $p$-adic higher residue pairing, compatible with the Hodge–Tate and crystalline filtrations.

This completes the bridge between the higher residue pairing on crystalline cohomology and the $p$-adic local system side. In particular, $(\mathbb{L}^i, \langle\ ,\ \rangle_{\mathbb{L}^i})$ provides a natural candidate for a $p$-adic version of Saito's polarized primitive form, suggesting a framework for $p$-adic TERP structures and $p$-adic Frobenius manifolds.